\documentclass[11pt,reqno]{amsart} 
\usepackage{amsmath,amssymb,rawfonts}
\usepackage{pb-diagram}
\newcommand{\impl}{\longmapsto} 
\newcommand{\ls}{\leqslant}
\newcommand{\be}{\begin{equation}}
\newcommand{\ee}{\end{equation}}
\newcommand{\cx}{c_{\text{\tiny{$X$}}}}
\newcommand{\cy}{c_{\text{\tiny{$Y$}}}}
\newcommand{\cz}{c_{\text{\tiny{$Z$}}}}

\newtheorem{theorem}{Theorem}[section]

\newtheorem{defi}[theorem]{Definition}
\newtheorem{prop}[theorem]{Proposition}
\newtheorem{coro}[theorem]{Corollary}

\newtheorem{ex}[theorem]{Example}

\begin{document} 
\title{Fixed and variable-basis fuzzy closure operators}
\author{Joaqu\'in Luna-Torres}
\dedicatory{Escuela de Matem\'aticas, Universidad Sergio Arboleda, Colombia}

\email{jluna@ima.usergioarboleda.edu.co}
\begin{abstract}
Closure operators are very useful tools in several areas of classical mathematics and in general  category theory. In fuzzy set  theory,  fuzzy closure operators have been studied by G. Gerla (1966). These works generally define a fuzzy subset as  a mapping from a  set $X$ to the real unit interval, as a complete and complemented lattice. More recently, Y. C. Kim (2003), F. G. Shi (2009), J. Fang and Y. Yue (2010)  propose theories of fuzzy closure systems and fuzzy closure operators in a more general settings, but still using complemented lattices.

The aim of this paper is to propose a more general theory of fixed and variable-basis fuzzy  closure operators, employing both categorical tools and the lattice theoretical fundations investigated by S. E. Rodabaugh (1999), where the lattices are usually non-complemented.  Besides, we construct topological categories in both cases.

\end{abstract}
\maketitle 
\baselineskip=1.7\baselineskip
\section*{0. Introduction}
It is well-known that the associated closure and interior operators provide equivalent descriptions of set-theoretic  topology; but this is not generally true in other categories, consequently it makes sense to define and study the notion of closure operators $C$ in the context of fuzzy set theory, where we can find categories in a lattice-theoretical context. 

Closure operators are very useful tools in  several areas of classical mathematics, and particularly in category theory. In fuzzy set  theory,  fuzzy closure operators have been studied by Gerla and others, (see e.g. \cite{GG}). These works generally define a fuzzy subset as  a mapping from a  set $X$ to the real unit interval $[0, 1]$, as a complete and complemented lattice. 

More recently, \cite{FY}, \cite{YCK} and \cite{FS}  propose  theories of fuzzy closure systems and fuzzy closure operators in a more general settings, but using complemented lattices.

The aim of this paper is to propose a more general theory of fixed and variable-basis fuzzy  closure operators, employing both categorical tools and the lattice theoretical fundations investigated in  \cite{SER} and \cite{HS}, where the lattices are usually non-complemented.

The paper is organized as follows: Following \cite{SER} and \cite{HS}  we introduce, in section $1$, the basic lattice theoretical fundations. In section $2$, we present the concept of fixed-basis fuzzy closure operators and then we construct a topological category $(\text{FBCO-SET},U)$, next  in section $3$ we present the concept of variable-basis fuzzy closure operators and then we construct a topological category $(\text{VBCO-SET},U)$. \\
In section 4, we study some additional stability properties of closure operators, that is idempotent and additive closure operators. Finally in section 5, we present some examples of various classes of closure maps. 

\section{From Lattice Theoretic Foundations}
Let $(L, \leq)$ be a complete, infinitely distributive lattice, i.e. $(L, \leq)$ is a partially ordered set such that for every subset $A\subset L$ the join $\bigvee A$ and the meet $\bigwedge A$ are defined, moreover $(\bigvee A) \wedge \alpha = \bigvee \{ a\wedge \alpha) \mid a \in A \}$ and\linebreak $(\bigwedge A) \vee \alpha = \bigwedge \{a\vee \alpha) \mid a \in A \}$
for every $\alpha \in L$.  In particular, $\bigvee L = \top$ and $\bigwedge L = \bot$ are respectively the universal upper and the universal lower bounds in $L$.
We assume that $\bot \ne \top$, i.e.\ $L$ has at least two elements.

\subsection{Complete quasi-monoidal lattices} 
The definition of complete quasi-monoidal lattices introduced by  S. E. Rodabaugh in \cite{SER} is the following:

A $cqm-$lattice (short for complete quasi-monoidal lattice) is a
triple\linebreak $(L,\leqslant,\otimes)$ provided with the
following properties
\begin{enumerate}
\item[(1)] $(L,\leqslant)$ is a complete lattice with upper bound $\top$ and lower bound $\bot$.
\item[(2)] $\otimes: L\times L  \rightarrow L$ is a binary operation satisfying the following axioms:
\begin{enumerate}
\item $\otimes$ is isotone in both arguments, i.e. $\alpha_1 \leqslant\alpha_2,\,\ \beta_1\leqslant \beta_2$ implies
$\alpha_1\otimes\beta_1\leqslant \alpha_2\otimes\beta_2$; \item  $\top$ is idempotent, i.e. $\top \otimes\top = \top$.
\end{enumerate}
\end{enumerate}
The category $CQML$ comprises the following data:
\begin{enumerate}
\item[(a)] {\bf Objects}: Complete quasi-monoidal lattices.
\item[(b)] {\bf Morphisms}: All $SET$ morphisms, between the above objects, which preserve $\otimes$ and $\top$ and arbitrary $\bigvee$.
\item[(c)] Composition and identities are taken from $SET$. 
\end{enumerate}
The category $LOQML$ is the dual of $CQML$, i.e. $LOQML=CQML^{op}$.
\subsection{$GL-$monoids}
A $GL-$monoid (see \cite{HS}) is a complete
lattice enriched with a further binary operation $\otimes$, i.e.\ a triple $(\L, \leq, \otimes)$ such that:
\begin{enumerate}
\item[(1)]
$\otimes$ is isotone, commutative and associative; 
\item[(2)]
$(\L,\leq,\otimes)$ is integral, i.e.\ $\top$ acts as the unity: $\alpha \otimes \top = \alpha$, $\forall \alpha \in \L$;
\item[(3)]
$\bot$ acts as the zero element in $(\L, \leq, \otimes)$, i.e.\ $\alpha\otimes \bot = \bot$, $\forall \alpha \in \L$;
\item[(4)]
$\otimes$ is distributive over arbitrary joins, i.e.\ $\alpha \otimes (\bigvee_{\lambda} \beta_{\lambda}) = \bigvee_{\lambda} (\alpha \otimes \beta_{\lambda})$,
$\forall \alpha \in \L, \forall \{ \beta_{\lambda} : \lambda \in I\} \subset \L$;
\item[(5)]
$(\L, \leq, \otimes)$ is divisible, i.e.\ $\alpha \leq \beta$ implies the existence of $\gamma \in \L$ such that $\alpha = \beta \otimes \gamma$.
\end{enumerate}
It is well known that every $GL-$monoid is residuated, i.e.\ there exists a further binary operation ``$\impl$'' (implication) on $\L$ satisfying the following condition:
$$\alpha \otimes \beta \leq \gamma \Longleftrightarrow \alpha \leq (\beta \impl \gamma) \qquad \forall \alpha, \beta, \gamma \in \L.$$
Explicitly the implication is given by
\[
\alpha \impl \beta = \bigvee \{ \lambda \in \L \mid \alpha \otimes \lambda \leq \beta \}.
\]
If $X$ is a set and $L$ is a $GL$-monoid (or a complete quasi-monoidal lattice), then the fuzzy powerset $L^X$ in an obvious way can be pointwise endowed with a structure of a $GL$-monoid (or of a complete quasi-monoidal lattice). In particular the $L$-sets $1_X$ and $0_X$ defined by $1_X (x)= \top$ and $0_X (x) = \bot$ $\forall x \in X$ are respectively the universal  upper and lower bounds in $L^X$.

\subsection{Powerset operator foundations}

We give the powerset operators, developed and justified in detail by S.E. Rodabaugh in \cite{SER} and \cite{SER1}.\
Let $f\in SET(X,Y)$,\,\ $L,M\in |CQML|$, \,\ $\phi\in LOQML(L,M)$, and $\wp(X)$,\ $\wp(Y),\ L^X,\ M^Y$ be the classical powerset of $X$, the classical powerset of $Y$, the $L$-powerset of $X$, and the $M$-powerset of $Y$, respectively. Then the following powerset operators are defined:
\begin{enumerate}
\item
$
f^{\rightarrow}:\wp(X)\rightarrow \wp(Y)\,\ \text{by}\,\ f^{\rightarrow}(A)=\{f(x)\mid x\in A\}
$
\item
$
f^{\leftarrow}:\wp(Y)\rightarrow \wp(X)\,\ \text{by}\,\ f^{\leftarrow}(B)=\{x\in X\mid f(x)\in B\}
$
\item
$
f_L^{\rightarrow}:L^X\rightarrow L^Y\,\ \text{by}\,\ f_L^{\rightarrow}(a)(y)=\bigvee_{f(x)=y} a(x)
$
\item
$
f_L^{\leftarrow}:L^Y\rightarrow L^X\,\ \text{by}\,\ f_L^{\leftarrow}(b)=b\circ f
$
\item
$
^{*}\phi:L\rightarrow M \,\ \text{by}\,\ ^{*}\phi(\alpha)=\bigwedge \{\beta\in M\mid \alpha\leqslant \phi^{op}(\beta)\}
$
\item
$
\langle^{*}\phi\rangle:L^X\rightarrow M^X \,\ \text{by}\,\ \langle^{*}\phi\rangle(a)={^{*}\phi}\circ a
$
\item
$
\langle\phi^{op}\rangle: M^X\rightarrow L^X\,\ \text{by}\,\ {\langle\phi^{op}\rangle}(b)=\phi^{op}\circ b
$
\item
$
\left( f,\Phi \right)^{\rightarrow}:L^X\rightarrow M^Y \,\ \text{by}\,\ \left( f,\Phi \right)^{\rightarrow}(a) = \bigwedge \{b\mid f_L^{\rightarrow}(a)\leqslant \left( \Phi^{op} \right)(b)\},
$
\item
$
\left( f,\Phi \right)^{\leftarrow} : M^Y\rightarrow L^X \,\ \text{by}\,\ \left( f,\Phi \right)^{\leftarrow}(b) =  \Phi^{op} \circ b\circ f,
$
in other words, that diagram
\[
\begin{diagram}
\node{X}\arrow{e,t}{f}\arrow{s,l}{\left( f,\Phi \right)^{\leftarrow}(b)}\node{Y}\arrow{s,r}{b}\\
\node{L}\node{M}\arrow{w,b}{\Phi^{op}}
\end{diagram}
\]
is commutative.
\end{enumerate}
Note that these operators were defined taking into account the Adjoint functor theorem. Consequently, we have that $f^{\rightarrow}$, \, $f_L^{\rightarrow}$, and $\left( f,\Phi \right)^{\rightarrow}$ are left adjoints of $f^{\leftarrow}$,\, $f_L^{\leftarrow}$, and $\left( f,\Phi \right)^{\leftarrow}$, respectively.
\section{Basic properties of fixed-basis fuzzy closure operators}
\begin{defi}\label{CO}
Given $L$ an object of the category $CQML$, a  closure operator $C$ of the category $SET$ (of  sets and functions between them) with respect to $L$ is given by a family $C=(c_{\text{\tiny{$X$}}})_{\text{\tiny{$X\in |SET|$}}}$ of maps $c_{\text{\tiny{$X$}}}:L^X\rightarrow L^X$ such that for every set $X$:
\begin{enumerate}
\item[($C_1$)] (Extension) $u\ls c_{\text{\tiny{$X$}}}(u)$\ for all $u\in L^X$;
\item[($C_2$)] (Monotonicity)  if $u\ls v$ in $L^X$, then $\cx(u)\ls \cx(v)$;
\item[($C_3$)] (Lower bound) $\cx(0_{\text{\tiny{$X$}}})=0_{\text{\tiny{$X$}}}$.
\end{enumerate}
\end{defi}
\begin{defi}
A fuzzy $C$-space is a pair $(X,\cx)$ where $X$ is a set and $\cx$ is a closure map on $X$.
\end{defi}
\begin{defi}
A function $f:X\rightarrow Y$ in $SET$ is said to be $C$-continuous if
\be\label{cont1}
f_L^{\rightarrow}\big( \cx(u) \big)\ls  c_{\text{\tiny{$Y$}}}\big(f_L^{\rightarrow}(u)\big) \,\ \text{for all}\,\ u\in L^X.
\ee
\end{defi}
\begin{prop}
The $C$-continuity condition can  equivalently be expressed as
\be\label{cont2}
\cx\big( f_L^{\leftarrow}(v)\big)\ls f_L^{\leftarrow}\big( c_{\text{\tiny{$Y$}}}(v)\big)\,\ \text{for all}\,\ v\in L^Y.
\ee
\end{prop}
\begin{proof}
Condition (\ref{cont1}) gives that, for $v\in L^Y$ and $u=f_L^{\leftarrow}(v)$,
\[
 f_L^{\rightarrow}\big(\cx (f_L^{\leftarrow}(v))\big)\ls c_{\text{\tiny{$Y$}}}\big(f_L^{\rightarrow}(f_L^{\leftarrow}(v))\big)\ls c_{\text{\tiny{$Y$}}}(v),
\]
then
\[
\cx\big( f_L^{\leftarrow}(v)\big)\ls f_L^{\leftarrow}\big( c_{\text{\tiny{$Y$}}}(v)\big).
\]
On the other hand, condition (\ref{cont2}) gives that, for $u\in L^X$ and $v=f_L^{\rightarrow}(u)$,
\[
\cx(u)\ls \cx\big(f_L^{\leftarrow}(f_L^{\rightarrow}(u))  \big)\ls f_L^{\leftarrow}\big( c_{\text{\tiny{$Y$}}}(f_L^{\rightarrow}(u))) \big),
\]
then
\[
f_L^{\rightarrow}\big( \cx(u) \big)\ls  c_{\text{\tiny{$Y$}}}\big(f_L^{\rightarrow}(u)\big).
\]
\end{proof}
\begin{prop}
Let $f:X\rightarrow Y$ and $g:Y\rightarrow Z$ be two $C$-continuous functions then the function $g\circ f$ is $C$-continuous.
\begin{proof}
Since $f:X\rightarrow Y$ is $C$-continuous, we have 
\[
f_L^{\rightarrow}\big( c_{\text{\tiny{$X$}}}(u) \big)\ls c_{\text{\tiny{$Y$}}}\big(f_L^{\rightarrow}(u)  \big)\,\ \text{for all}\, u\in L^X,
\]
it follows that
\[
g_L^{\rightarrow}\Big(  f_L^{\rightarrow}\big( \cx(u) \big)\Big)\ls g_L^{\rightarrow}\Big( \cy\big(f_L^{\rightarrow}(u)  \big)\Big),
\]
now, by the $C$-continuity  of $g$,
\[g_L^{\rightarrow}\Big(\cy\big(  f_L^{\rightarrow}(u) \big)\Big)\ls \cz\Big( g_L^{\rightarrow}\big(f_L^{\rightarrow}(u)  \big)\Big)
\]
therefore 
\[
\big(g\circ f\big)_L^{\rightarrow}\big( \cx(u) \big)\ls \cz \big(g\circ f\big)_L^{\rightarrow}(u)  
\]
\end{proof}
\end{prop}
as a consequence we obtain

\begin{defi}
The category FBCO-SET of $C$-spaces comprises de following data:
\begin{enumerate}
\item[(1)] {\bf Objects}: pairs $(X,\cx)$, where $X$ is a set and $\cx$ is a clousure map on $X$.
\item[(2)]{\bf Morphisms}: Functions which are $C$-continuous.
\end{enumerate}
\end{defi}
\subsection{ The lattice structure of all closure operators.}
We consider the collection
\[C\big(SET,L\big)
\]
of all closure operators on $SET$ with respect to the complete quasi-monoidal lattice $L$. It is ordered by
$$C \ls D\Leftrightarrow \cx (u) \ls d_{\text{\tiny{$X$}}}(u),\, \text{for all set $X$, and for all}\, u\in L^X.$$
This way $C\big(SET,L\big)$ inherents a lattice structure from $L$:

\begin{prop}
Every family $\big(C_{\text{\tiny{$\lambda$}}}\big)_{\text{\tiny{$\lambda \in \Lambda$}}}$ in $C\big(SET,L\big)$ has a join $\bigvee \limits_{\lambda \in \Lambda} C_{\text{\tiny{$\lambda$}}}$ and a meet $\bigwedge\limits_{\lambda \in \Lambda} C_{\text{\tiny{$\lambda$}}}$ in $C\big(SET,L\big)$.
The discrete closure operator
\[ C_{\text{\tiny{$D$}}}= \big(c_{\text{\tiny{$D_X$}}}\big)_{X \in |SET|}
\]
is the least element in $C\big(SET,L\big)$, and the trivial closure operator
\[ C_{\text{\tiny{$T$}}}= \big(c_{\text{\tiny{$T_X$}}}\big)_{X \in |SET|}\qquad \text{with}\qquad \big(c_{\text{\tiny{$T_X$}}}\big)(u)=\begin{cases}&1_{\text{\tiny{$X$}}} \; \text{for all} \,\,\ u\neq 0 \\ & 0_{\text{\tiny{$X$}}}\; \text{if} \,\,\ u=0_{\text{\tiny{$X$}}}
\end{cases}
\]
is the largest one.

\end{prop}

\begin{proof}
For $\Lambda \neq \emptyset $, let $\tilde{C}=\bigwedge\limits_{\lambda \in \Lambda} C_{\text{\tiny{$\lambda$}}}$, then  
$$\tilde{C}_{\text{\tiny{$X$}}}=\bigwedge\limits_{\lambda \in \Lambda} C_{\text{\tiny{$\lambda_{X}$}}},$$ where $X$ is an arbitrary set, satisfies 
\begin{itemize}
\item $u \ls \tilde{c}_{\text{\tiny{$X$}}}(u)$, because $u\ls c_{\text{\tiny{$\lambda_{X}$}}}(u)$ for all $u \in L^X$ and for all $\lambda \in \Lambda.$ 
\item If $u_1 \ls u_2$ in $L^X$ then $ c_{\text{\tiny{$\lambda_{X}$}}}(u_1)\ls c_{\text{\tiny{$\lambda_{X}$}}}(u_2)$ for all $\lambda \in \Lambda$, therefore $\tilde{c}_{\text{\tiny{$X$}}}(u_1)\ls \tilde{c}_{\text{\tiny{$X$}}}(u_2).$
\item  Since $c_{\text{\tiny{$\lambda_{X}$}}}(0_X)=0_X$ for all $\lambda \in \Lambda$, we have that $\tilde{c}_{\text{\tiny{$X$}}}(0_X)=0_X.$

\end{itemize}

Similary $\bigvee\limits_{\lambda \in \Lambda} C_{\text{\tiny{$\lambda_{X}$}}}$, $C_{\text{\tiny{$T_X$}}}$ and $C_{\text{\tiny{$D_X$}}}$ are closure operators.

\end{proof}
Consequently,
\begin{coro}\label{cl}
For every set $X$
$$CL(X)=\{c_{\text{\tiny{$X$}}} \mid  c_{\text{\tiny{$X$}}} \, \text{is a closure map on X}\}$$ is a complete lattice.

\end{coro}
\subsection{Initial closure operators}
Let $(Y,c_{\text{\tiny{$Y$}}})$ be an object of the category FBCO-SET, and let $X$ be a set. For each function $f:X\rightarrow Y$ we define on $X$ the map 

\be \label{cli}
 c_{\text{\tiny{$X_f$}}}:L^X\rightarrow L^X \qquad  \text{  by } \qquad c_{\text{\tiny{$X_f$}}}=f_{\text{\tiny{$L$}}}^{\text{\tiny{$\leftarrow$}}} \circ c_{\text{\tiny{$Y$}}}\circ f_{\text{\tiny{$L$}}}^{\text{\tiny{$\rightarrow$}}}
\ee
i.e. the following diagram is conmutative 
\[
\begin{diagram}
\node{L^X}\arrow{e,t}{f_{\text{\tiny{$L$}}}^{\text{\tiny{$\rightarrow$}}}} \arrow{s,..}{ c_{\text{\tiny{$X_f$}}}}\node{L^Y}\arrow{s,r}{c_{\text{\tiny{$Y$}}}}\\
\node{L^X}\node{L^Y}\arrow{w,b}{f_{\text{\tiny{$L$}}}^{\text{\tiny{$\leftarrow$}}}}
\end{diagram}
\]

\begin{prop}\label{ci}
The map (\ref{cli}) is a closure map on $X$ for which the function $f$ is C-continuous. 
\end{prop}
\begin{proof}\
\begin{enumerate}
\item[1)](Extension) $c_{\text{\tiny{$X_f$}}}(u)= f_{\text{\tiny{$L$}}}^{\text{\tiny{$\leftarrow$}}}\Big(c_{\text{\tiny{$Y$}}}\big(f_{\text{\tiny{$L$}}}^{\text{\tiny{$\rightarrow$}}}(u)\big) \Big)\geqslant f_{\text{\tiny{$L$}}}^{\text{\tiny{$\leftarrow$}}}\big(f_{\text{\tiny{$L$}}}^{\text{\tiny{$\rightarrow$}}}(u)\big)\geqslant u$,\\  therefore $u\ls c_{\text{\tiny{$X_f$}}}(u)$
\item[2)](Monotonicity) $u_1\ls u_2$ in $L^X$ implies $f_{\text{\tiny{$L$}}}^{\text{\tiny{$\rightarrow$}}}(u_1)\ls f_{\text{\tiny{$L$}}}^{\text{\tiny{$\rightarrow$}}}(u_2)$\\
 then
 $c_{\text{\tiny{$Y$}}}\big(f_{\text{\tiny{$L$}}}^{\text{\tiny{$\rightarrow$}}}(u_1)\big)\ls c_{\text{\tiny{$Y$}}}\big(f_{\text{\tiny{$L$}}}^{\text{\tiny{$\rightarrow$}}}(u_2)\big)$, consequently $$c_{\text{\tiny{$X_f$}}}(u_1)=\Big(f_{\text{\tiny{$L$}}}^{\text{\tiny{$\leftarrow$}}}\Big(c_{\text{\tiny{$Y$}}} \big(f_{\text{\tiny{$L$}}}^{\text{\tiny{$\rightarrow$}}}(u_1)\big)\Big)\Big)\ls \Big(f_{\text{\tiny{$L$}}}^{\text{\tiny{$\rightarrow$}}}\Big(c_{\text{\tiny{$Y$}}} \big(f_{\text{\tiny{$L$}}}^{\text{\tiny{$\rightarrow$}}}(u_2)\big)\Big)\Big)=c_{\text{\tiny{$X_f$}}}(u_2).$$
 \item[3)](Lower bound)  $c_{\text{\tiny{$X_f$}}}(0_X)= f_{\text{\tiny{$L$}}}^{\text{\tiny{$\leftarrow$}}}\Big(c_{\text{\tiny{$Y$}}}\big(f_{\text{\tiny{$L$}}}^{\text{\tiny{$\rightarrow$}}}(0_X)\big) \Big)=0_X$.\\
 Finally $$f_{\text{\tiny{$L$}}}^{\text{\tiny{$\rightarrow$}}}\Big(c_{\text{\tiny{$X_f$}}}(u)\Big)=f_{\text{\tiny{$L$}}}^{\text{\tiny{$\rightarrow$}}}\Big(f_{\text{\tiny{$L$}}}^{\text{\tiny{$\leftarrow$}}}\Big(c_{\text{\tiny{$Y$}}}\big(f_{\text{\tiny{$L$}}}^{\text{\tiny{$\rightarrow$}}}(u)\big) \Big)\Big)\ls c_{\text{\tiny{$Y$}}}\Big(f_{\text{\tiny{$L$}}}^{\text{\tiny{$\rightarrow$}}}(u)\Big) \,\text{for all}\, u\in L^X.$$
 \end{enumerate}
\end{proof}
It is clear that $c_{\text{\tiny{$X_f$}}}$ is the finner map on $L^X$ for which the function $f$ is C-continuous, more precisaly.

\begin{prop}\label{unique}
Let $(Z,c_{\text{\tiny{$Z$}}})$ and $(Y,c_{\text{\tiny{$Y$}}})$ be objects of FBCO-SET, and let $X$ be a set. For each function  $g:Z\rightarrow X$  and for  $f:(X,c_{\text{\tiny{$X_{f}$}}})\rightarrow (Y,c_{\text{\tiny{$Y$}}})$ a $C$-continuous function, $g$  is $C$-continuous if and only if $f \circ g$ is $C$-continuous.
\end{prop}

\begin{proof}
Suppose that $f\circ g$ is $C$-continuous, i. e.
$$c_{\text{\tiny{$Z$}}}\big((f\circ g)_{\text{\tiny{$L$}}}^{\text{\tiny{$\leftarrow$}}}(v)\big)\leqslant \big( (f\circ g)_{\text{\tiny{$L$}}}^{\text{\tiny{$\leftarrow$}}} c_{\text{\tiny{$Y$}}}(v) \big)$$
 for all $v\in L^Y$. Then, for all $u\in L^X$, we have
\begin{align*}
 g_{\text{\tiny{$L$}}}^{\text{\tiny{$\leftarrow$}}}\big(c_{\text{\tiny{$X_{f}$}}}(u)\big)&=g_{\text{\tiny{$L$}}}^{\text{\tiny{$\leftarrow$}}}(\big( f_{\text{\tiny{$L$}}}^{\text{\tiny{$\leftarrow$}}}\circ c_{\text{\tiny{$Y$}}}\circ f_{\text{\tiny{$L$}}}^{\text{\tiny{$\rightarrow$}}})(u)\big)= (f\circ g)_{\text{\tiny{$L$}}}^{\text{\tiny{$\leftarrow$}}}\big( c_{\text{\tiny{$Y$}}}( (f_{\text{\tiny{$L$}}}^{\text{\tiny{$\rightarrow$}}}(u)) \big)\\
 &\geqslant c_{\text{\tiny{$Z$}}}\big( (f\circ g)_{\text{\tiny{$L$}}}^{\text{\tiny{$\leftarrow$}}}f_{\text{\tiny{$L$}}}^{\text{\tiny{$\rightarrow$}}}(u) ) \big)= c_{\text{\tiny{$Z$}}}\big( (g_{\text{\tiny{$L$}}}^{\text{\tiny{$\leftarrow$}}}\circ f_{\text{\tiny{$L$}}}^{\text{\tiny{$\leftarrow$}}}\circ f_{\text{\tiny{$L$}}}^{\text{\tiny{$\rightarrow$}}})(u)\big)\geqslant c_{\text{\tiny{$Z$}}}\big( g_{\text{\tiny{$L$}}}^{\text{\tiny{$\leftarrow$}}}(u) \big).\\
\end{align*}
\end{proof}
\begin{prop}\label{ss}
Let $X$ be a set, let $(Y_j,c_{\text{\tiny{$Y_j$}}})$ be a family of fuzzy $C$-spaces, where $j\in J$ for some indexed set $J$, and let $f_j:X\rightarrow Y_j$ be functions. Then the structured source $(X,f_j:X\rightarrow Y_j)$  w.r.t the forgetful functor $U$ from FBCO-SET to SET has a unique initial lift $\big( (X, \hat{c}_{\text{\tiny{$X$}}})\rightarrow (Y_j,c_{\text{\tiny{$Y_j$}}})\big)$, where $\hat{c}_{\text{\tiny{$X$}}}$ is the join  $\bigvee \limits_{\lambda \in \Lambda} c_{\text{\tiny{$f_j$}}}$ of all initial closure maps $c_{\text{\tiny{$f_j$}}}$ w.r.t. $f_j$, where $j\in J$.
\end{prop}
\begin{proof}
We must show that for every object $(Z, c_{\text{\tiny{$Z$}}})$  of FBCO-SET,  each function $g:Z\rightarrow X$ is C-continuous iff $f_j\circ g$ is C-continuous, for all $j\in J$. \\
In fact,
\begin{align*}
g_{\text{\tiny{$L$}}}^{\text{\tiny{$\leftarrow$}}}\Big(\bigvee \limits_{j \in J} c_{\text{\tiny{$f_j$}}}(u)\Big)&= g_{\text{\tiny{$L$}}}^{\text{\tiny{$\leftarrow$}}}\Big(\bigvee \limits_{j \in J}\big(f_{j}^{\text{\tiny{$\leftarrow$}}}\circ c_{\text{\tiny{$f_j$}}}\circ f_{j}^{\text{\tiny{$\rightarrow$}}}\big)(u)\Big)\\
&=\bigvee \limits_{j \in J}g_{\text{\tiny{$L$}}}^{\text{\tiny{$\leftarrow$}}}\Big(\big(f_{j}^{\text{\tiny{$\leftarrow$}}}\circ c_{\text{\tiny{$f_j$}}}\circ f_{j}^{\text{\tiny{$\rightarrow$}}}\big)(u)\Big)\\
&\geqslant \big( f_{j}\circ g\big)^{\text{\tiny{$\leftarrow$}}} \big( c_{\text{\tiny{$Y_j$}}} (f_{j}^{\text{\tiny{$\rightarrow$}}}(u)) \big) \\
&\geqslant c_{\text{\tiny{$Z$}}}\Big( ( f_{j}\circ g\big)^{\text{\tiny{$\leftarrow$}}}\big( f_{j}^{\text{\tiny{$\rightarrow$}}}(u) \big) \Big)\\
&\geqslant c_{\text{\tiny{$Z$}}}\big(g_{\text{\tiny{$L$}}}^{\text{\tiny{$\leftarrow$}}}(u) \big).
\end{align*}
\end{proof}
Now, remember that 
\begin{defi} $\left( \text{Ad\'amek \cite{AHS} and Rodabaugh \cite{SER1}}\right)$  Category $\mathcal A$ is topological with regard to $\mathcal X$ and  functor $ V:\mathcal A\rightarrow \mathcal X$ iff  eah $V$-structured source in $\mathcal X$ has a unique, initial $V$-lift in $\mathcal A$. We may also say that $\mathcal A$ is topological over $\mathcal X$ with regard to functor $V$.
\end{defi}
   As a consequence of corollary (\ref{cl}), proposition (\ref{ci}) and proposition (\ref{ss}), we obtain
          
  \begin{theorem} \label{topcat}
  The concrete category $(\text{FBCO-SET},\mathcal O)$ over $SET$ is a topological category.
  \end{theorem}              
                      
\section{Basic properties of variable-basis closure operators}
In this section we consider a subcategory $\mathcal{D}$ of CQML in order to construct fuzzy variable-basis closure operators on the category $SET\times\mathcal{D}$ that has as objects all pairs $(X,L)$, where $X$ is a set and $L$ is an object of $\mathcal{D}$, as morphisms from $(X, L)$ to $(Y,M)$ all pairs of maps $(f, \phi)$ with $f\in SET(X,L)$ and $ \phi \in CQML(L,M)$, identities given by $id_{\text{\tiny{$(X,L)$}}} =(id_{\text{\tiny{$X$}}} ,id_{\text{\tiny{$L$}}}) $, and composition defined by

\[
(f, \phi ) \circ(g, \psi)=(f\circ g, \phi\circ \psi).
\]

\begin{defi} A closure operator of the category $SET\times\mathcal{D}$ is given by a family

$C=(c_{\text{\tiny{$XL$}}})_{\text{\tiny{$(X,L)$}} \in \left| SET\times\mathcal{D} \right|}$ of maps 
$c_{\text{\tiny{$XL$}}}:L^{\text{\tiny{$X$}}}\longrightarrow L^{\text{\tiny{$X$}}}$ that satisfies the requeriment:

\begin{enumerate}
\item[($C_1$)] (Extension) $u\ls c_{\text{\tiny{$X$}}}(u)$\ for all $u\in L^X$;
\item[($C_2$)] (Monotonicity)  if $u\ls v$ in $L^X$, then $\cx(u)\ls \cx(v)$;
\item[($C_3$)] (Lower bound) $\cx(0_{\text{\tiny{$X$}}})=0_{\text{\tiny{$X$}}}$.
\end{enumerate}
\end{defi}
                         
 \begin{defi} A fuzzy variable-basis $C$-space is a triple $(X,L, c_{\text{\tiny{$XL$}}})$, where $(X,L)$ is an
object of $SET\times\mathcal{D}$ and $c_{\text{\tiny{$XL$}}}$ is a closure map on $(X,L)$.
\end{defi}

\begin{defi}
 A morphism $(f, \phi):(X,L) \longrightarrow (Y,M)$ in $SET\times\mathcal{D}$ is said to be fuzzy $c$-continuous if

\begin{equation}\label{ec1}
(f, \phi)^{\rightarrow}\big(c_{\text{\tiny{$XL$}}}(u)\big) \ls c_{\text{\tiny{$YM$}}}\Big((f,\phi)^\rightarrow (u)\Big)\; \text{for all} \; u \in L^{\text{\tiny{$X$}}} . 
\end{equation}
\end{defi}

\begin{prop}
Condition (\ref{ec1}) is equivalent to 
\begin{equation}\label{ec2}
c_{\text{\tiny{$XL$}}}\Big((f,\phi)^\leftarrow (v)\Big)\ls (f,\phi)^\leftarrow\big(c_{\text{\tiny{$YM$}}}(v) \big) \quad \text{for all} \quad v \in M^{\text{\tiny{$Y$}}} .
\end{equation}
\end{prop}
\begin{proof}
Condition (\ref{ec1}) gives that, for $v \in M^{\text{\tiny{$Y$}}}$ and $u=(f,\phi)^\leftarrow (v)$, 

\begin{align*}
(f,\phi)^{\rightarrow}\big(c_{\text{\tiny{$XL$}}}(u)\big)&=(f,\phi)^{\rightarrow}\big(c_{\text{\tiny{$XL$}}}(f,\phi)^\leftarrow (v)\big)\\
&\ls c_{\text{\tiny{$YM$}}}\Big((f,\phi)^{\rightarrow}(f,\phi)^\leftarrow (v)\Big) \ls c_{\text{\tiny{$YM$}}}(v),
\end{align*}
then
\[
c_{\text{\tiny{$XL$}}}\big((f,\phi)^\leftarrow (v)\big)\ls (f,\phi)^\leftarrow \big(c_{\text{\tiny{$YM$}}}(v)\big).
\]
On the other hand, condition\,\  (\ref{ec2}) gives that, for $ u \in L^{\text{\tiny{$X$}}}$ 

\begin{align*}
c_{\text{\tiny{$XL$}}}(u)&\ls c_{\text{\tiny{$XL$}}}\Big((f,\phi)^\leftarrow (f,\phi)^\rightarrow (u)\Big)\\
&\ls (f,\phi)^\leftarrow \Big(c_{\text{\tiny{$YM$}}}\big((f,\phi)^\rightarrow (u)\big)\Big),
\end{align*}
then $(f,\phi)^\rightarrow(c_{\text{\tiny{$XL$}}}(u))\ls c_{\text{\tiny{$YM$}}}\big((f,\phi)^\rightarrow (u)\big)$
\end{proof}

\begin{prop}
Consider two fuzzy $c$-continuous morphisms \linebreak $(f, \phi):(X,L) \longrightarrow (Y,M)$ and $(g,\psi):(Y,M)\longrightarrow (Z,N)$ be , then the morphism $(g, \psi)\circ (f,\phi)$ is fuzzy $c$-continuous.
\end{prop}

\begin{proof}
Since $(f, \phi):(X,L) \longrightarrow (Y,M)$ is $c$-continuous we have
\[(f,\phi)^\rightarrow(c_{\text{\tiny{$XL$}}}(u))\ls c_{\text{\tiny{$YM$}}}\big((f,\phi)^\rightarrow (u)\big) \quad \text{for all} \quad u \in L^{\text{\tiny{$X$}}} ,
\]
 it follows that
$$(g, \psi)^\rightarrow \Big((f,\phi)^\rightarrow \big(c_{\text{\tiny{$XL$}}}(u)\big) \Big)\ls (g, \psi)^\rightarrow \big(c_{\text{\tiny{$YM$}}}\big((f,\phi)^\rightarrow (u)\big)\big),$$ 
now, by the fuzzy $c$-continuity of $(g,\psi)$,
$$(g, \psi)^\rightarrow \Big( c_{\text{\tiny{$YM$}}} \big((f,\phi)^\rightarrow (u)\big)  \Big)\ls c_{\text{\tiny{$ZN$}}}\Big((g, \psi)^\rightarrow \big((f,\phi)^\rightarrow (u)\big)\Big),$$ therefore 
\[
\Big((g, \psi)\circ (f,\phi)\Big)^\rightarrow\big(c_{\text{\tiny{$XL$}}}(u)\big)\ls c_{\text{\tiny{$ZN$}}}\Big(\big((g, \psi)\circ (f,\phi)\big)^\rightarrow (u)\Big).\]
\end{proof}
    
    As a consequence we obtain
 \begin{defi} The category VBCO-SET that has as objects all triples $(X,L,c_{\text{\tiny{$XL$}}})$ where $(X,L)$ is an object of $SET\times\mathcal{D}$ and $c_{\text{\tiny{$XL$}}}:L^{\text{\tiny{$X$}}}\longrightarrow L^{\text{\tiny{$X$}}}$ is a fuzzy closure map, as morphisms from $(X,L,c_{\text{\tiny{$XL$}}})$ to $(Y,M,c_{\text{\tiny{$YM$}}})$ all pairs of fuzzy $c$-continuous functions $(f,\phi):(X,L,c_{\text{\tiny{$XL$}}})\longrightarrow (Y,M,c_{\text{\tiny{$YM$}}})$, identities and composition as in $SET\times D$
 \end{defi} 
 
  \subsection{Initial variable-basis closure operator}   
 Let $(Y,M,c_{\text{\tiny{$YM$}}})$ be an object of the category   VBCO-SET and let $(X,L)$ be an object of the category $SET\times\mathcal{D}$.\\
  For each morphism $(f,\phi):(X,L)\longrightarrow(Y,M)$ in $SET\times\mathcal{D}$ we define on $(X,L)$ the map
 $\hat{c}_{\text{\tiny{$XL$}}}:L^{\text{\tiny{$X$}}}\longrightarrow L^{\text{\tiny{$X$}}}$ by 
 \begin{equation}\label{ec*}
 \hat{c}_{\text{\tiny{$XL$}}}=(f,\phi)^\leftarrow \circ c_{\text{\tiny{$YM$}}}\circ (f,\phi)^\rightarrow 
 \end{equation}
  i.e the following diagram is conmutative
 \[
\begin{diagram}
\node{L^{\text{\tiny{$X$}}}}\arrow{e,t}{(f,\phi)^\rightarrow}\arrow{s,..}{\hat{c}_{\text{\tiny{$XL$}}}}\node{M^{\text{\tiny{$Y$}}}}\arrow{s,r}{c_{\text{\tiny{$YM$}}}}\\
\node{L^{\text{\tiny{$X$}}}}\node{M^{\text{\tiny{$Y$}}}}\arrow{w,b}{(f,\phi)^\leftarrow}
\end{diagram}
\]                           
 \begin{prop}\label{prop.penult}
 The map (\ref{ec*}) is a closure map on $(X,L)$ for which the morphism $(f,\phi)$ is fuzzy $c$-continuous 
           \end{prop}
           
           \begin{proof}\
 \begin{enumerate}
 \item[1)] (Extension) For every $ u \in L^{\text{\tiny{$X$}}}$,
\[\hat{c}_{\text{\tiny{$XL$}}}(u)=(f,\phi)^\leftarrow\Big( c_{\text{\tiny{$YM$}}}\big( (f,\phi)^\rightarrow (u) \big)\Big)\geqslant (f,\phi)^\leftarrow (f,\phi)^\rightarrow (u)\geqslant u,
\]
 therefore $u\ls \hat{c}_{\text{\tiny{$XL$}}}(u).$
 
 \item[2)](Monotonicity) $u_1\ls u_2$ in $L^{\text{\tiny{$X$}}}$ implies $(f,\phi)^\rightarrow (u_1)\ls$  $(f,\phi)^\rightarrow (u_2),$
  then $c_{\text{\tiny{$YM$}}}\Big((f,\phi)^\rightarrow (u_1)\Big)\ls c_{\text{\tiny{$YM$}}}\Big((f,\phi)^\rightarrow (u_2)\Big),$\\
   consequentently\\
\begin{align*}
\hat{c}_{\text{\tiny{$XL$}}}(u_1)&=(f,\phi)^\leftarrow \Big(c_{\text{\tiny{$YM$}}}\big((f,\phi)^\rightarrow (u_1)\big)\Big)\\
 &\ls  (f,\phi)^\leftarrow \Big(c_{\text{\tiny{$YM$}}}\big((f,\phi)^\rightarrow (u_2)\big)\Big)=\hat{c}_{\text{\tiny{$XL$}}}(u_2)
\end{align*}
           
 \item[3)] (Lower bound) 
 \[ 
 \hat{c}_{\text{\tiny{$XL$}}}(0_X)=(f,\phi)^\leftarrow \Big(c_{\text{\tiny{$YM$}}}\big((f,\phi)^\rightarrow (0_X)\big)\Big)=0_X.
 \]
Finally, 
\begin{align*}
(f,\phi)^\rightarrow \big(\hat{c}_{\text{\tiny{$XL$}}}(u)\big)&=(f,\phi)^\rightarrow \Big((f,\phi)^\leftarrow \big(c_{\text{\tiny{$YM$}}}((f,\phi)^\rightarrow (u))\big)\Big)\\
&\ls c_{\text{\tiny{$YM$}}} \big( (f,\phi)^\rightarrow (u))\big),\, \text{ for all}\, u\in L^{\text{\tiny{$X$}}}.
\end{align*}         
 \end{enumerate}
           \end{proof}   
           It is clear that $\hat{c}_{\text{\tiny{$XL$}}}$ is the finest map on $L^{\text{\tiny{$X$}}}$ for which the morphism $(f,\phi)$ es fuzzy $c$-continuous, more precisely
           
           \begin{prop}\label{prop.ultima}
           Let $(Z,N,c_{\text{\tiny{$ZN$}}})$ and $(Y,M,c_{\text{\tiny{$YM$}}})$ be objects of the category VBCO-SET and let $(X,L)$ be an object of $SET\times\mathcal{D}$. For each morphism $(g,\psi):(Z,N)\longrightarrow (X,L)$ of $SET\times\mathcal{D}$ and for $(f,\phi):(X,\hat{c}_{\text{\tiny{$XL$}}})\longrightarrow (Y,c_{\text{\tiny{$YN$}}})$ a fuzzy $c$-continuous morphism, $(g,\psi)$ is a fuzzy $c$-continuous morphism if and only if $(f,\phi)\circ(g,\psi)$ is fuzzy $c$-continuous.
           \end{prop}

          \begin{proof}
Suppose that $(f,\phi)\circ(g,\psi)$ is fuzzy $c$-continuous, that is 
\[c_{\text{\tiny{$ZN$}}}\Big(\big((f,\phi)\circ(g,\psi)\big)^\leftarrow(v)\Big)\ls \big((f,\phi)\circ(g,\psi)\big)^\leftarrow \big(c_{\text{\tiny{$YN$}}}(v)\big)\, \text{for  all}\, v\in L^{\text{\tiny{$Y$}}}.
\]
Then for all $u\in L^X$,  we have 
        \begin{align*}
          (g,\psi)^\leftarrow \big(\hat{c}_{\text{\tiny{$XL$}}}(u)\big)&=(g,\psi)^\leftarrow \Big((f,\phi)^\leftarrow \big(c_{\text{\tiny{$YM$}}}((f,\phi)^\rightarrow (u))\big)\Big)\\
         &= \big((f,\phi)\circ (g,\psi)\big)^\leftarrow \big(c_{\text{\tiny{$YM$}}}((f,\phi)^\rightarrow (u))\big)\\
     &\geqslant c_{\text{\tiny{$ZN$}}} \big((f,\phi)\circ (g,\psi)\big)^\leftarrow \big((f,\phi)^\rightarrow (u)\big)\\
     & =c_{\text{\tiny{$ZN$}}}(g,\psi)^\leftarrow \big((f,\phi)^\leftarrow (f,\phi)^\rightarrow(u)\big)\\
     &\geqslant\big((g,\psi)^\leftarrow (u)\big)
      \end{align*}
          \end{proof}

 As in  theorem (\ref{topcat}), we have
          
          \begin{theorem} The concrete category $(\text{VBCO-SET}, U)$ over $SET\times\mathcal{D}$ is a topological category.
          \end{theorem}
          
          \subsection{Closed and dense fuzzy sets}
          \begin{defi}
          An $L$-fuzzy subset $u$ of $X$ is called $c$-closed in $(X,L)$ if it is equal to it is closure, i.e, $c_{\text{\tiny{$XL$}}}(u)=u.$ The fuzzy $c$-continuity condition (\ref{ec1}) implies that $c$-closedness is preserved by inverse images.
          \end{defi}
          \begin{prop} Let $(f,\phi):(X,L,c_{\text{\tiny{$XL$}}})\longrightarrow(Y,M,c_{\text{\tiny{$YM$}}})$ be a morphism in VBCO-SET. If $v \in M^{\text{\tiny{$Y$}}}$ is $c$-closed then $(f,\phi)^\leftarrow(v)$ is $c$-closed in $(X,L).$
          \end{prop}

   \begin{proof}
   If $v=c_{\text{\tiny{$YM$}}}(v),$ for $v\in M^{\text{\tiny{$Y$}}},$ then $(f,\phi)^\leftarrow(v)=(f,\phi)^\leftarrow\big(c_{\text{\tiny{$YM$}}}(v)\big)\geqslant c_{\text{\tiny{$XL$}}}\big((f,\phi)^\leftarrow(v)\big),$ so $c_{\text{\tiny{$XL$}}}\big((f,\phi)^\leftarrow(v)\big)=(f,\phi)^\leftarrow(v)$
   \end{proof}
   
   \begin{defi}
   An $L$-fuzzy subset $u$ of $X$ is called $c$-dense in $(X,L)$ if it is $c$-closure is $1_X.$
   \end{defi}

   \begin{prop}
   Let $(f,\phi):(X,L,c_{\text{\tiny{$XL$}}})\longrightarrow(Y,M,c_{\text{\tiny{$YM$}}}) $ be an epimorphism in VBCO-SET. If $u\in L^{\text{\tiny{$X$}}}$ is $(f,\phi)^\rightarrow (u)$ is $c$-dense in $(Y,M)$.
   \end{prop}   
   \begin{proof} It $c_{\text{\tiny{$XL$}}}(u)=1_X$ then 
   
   $1_Y=(f,\phi)^\rightarrow (1_X)= (f,\phi)^\rightarrow \big(c_{\text{\tiny{$XL$}}}(u)\big)\ls c_{\text{\tiny{$YM$}}}\big((f,\phi)^\rightarrow(u)\big)$
   \end{proof}  
 
    \subsection{$C$-closed morphisms}
    \begin{defi}
    A morphism $(f,\phi):(X,L,c_{\text{\tiny{$XL$}}})\longrightarrow(Y,M,c_{\text{\tiny{$YM$}}})$ between fuzzy variable-basis $c$-spaces is fuzzy $c$-closed if
    \begin{equation}\label{ec3}
  c_{\text{\tiny{$YM$}}}\big((f,\phi)^\rightarrow(u)\big)\ls(f,\phi)^\rightarrow\big( c_{\text{\tiny{$XL$}}}(u)\big)\quad \text{for all}\quad u \in L^{\text{\tiny{$X$}}}.
    \end{equation}
    \end{defi} 
    \begin{prop} Let $(f,\phi):(X,L,c_{\text{\tiny{$XL$}}})\longrightarrow(Y,M,c_{\text{\tiny{$YM$}}})$ and\linebreak $(g,\psi):(Y,M,c_{\text{\tiny{$YM$}}})\longrightarrow(Z,N,c_{\text{\tiny{$ZN$}}})$ be two fuzzy $c$-closed morphisms, then the morphism $(f,\phi)\circ (g,\psi)$ is fuzzy $c$-closed.
        \end{prop}
        
        \begin{proof} Since $(f,\phi):(X,L,c_{\text{\tiny{$XL$}}})\longrightarrow(Y,M,c_{\text{\tiny{$YM$}}})$ is fuzzy $c$-closed, we have $$c_{\text{\tiny{$YM$}}}\big((f,\phi)^\rightarrow(u)\big)\ls(f,\phi)^\rightarrow\big( c_{\text{\tiny{$XL$}}}(u)\big)\quad \text{for all}\quad u \in L^{\text{\tiny{$X$}}},$$
         it follows that 
         $$(g,\psi)^\rightarrow\Big( c_{\text{\tiny{$YM$}}}\big((f,\phi)^\rightarrow(u)\big)\Big)\ls (g,\psi)^\rightarrow \Big((f,\phi)^\rightarrow\big(c_{\text{\tiny{$XL$}}}(u)\big)\Big)$$
        
  now, by the fuzzy $c$-closedness of $(g,\psi),$
  $$c_{\text{\tiny{$ZN$}}}\big((g,\psi)^\rightarrow(v)\big)\ls(g,\psi)^\rightarrow\big(c_{\text{\tiny{$YM$}}}(v)\big)\quad \text{for all}\quad v \in L^{\text{\tiny{$Y$}}},$$ in particular for $v=(f,\phi)^\rightarrow(u),$ 
  $$c_{\text{\tiny{$ZN$}}}\big((g,\psi)^\rightarrow((f,\phi)^\rightarrow(u))\big)\ls(g,\psi)^\rightarrow\big(c_{\text{\tiny{$YM$}}}((f,\phi)^\rightarrow(u))\big)$$
  therefore
  $$c_{\text{\tiny{$ZN$}}}\Big(\big((g,\psi)\circ (f,\phi)\big)^\rightarrow(u)\Big)\ls (\big((g,\psi)\circ (f,\phi)\big)^\rightarrow \big(c_{\text{\tiny{$XL$}}}(u)\big)$$
       \end{proof}
If we replace in the category VBCO-SET fuzzy $c$-continuous morphisms by fuzzy $c$-closed morphisms, we obtain another topological category. 
The morphisms $(f,\phi):(X,L,c_{\text{\tiny{$XL$}}})\longrightarrow(Y,M,c_{\text{\tiny{$YM$}}})$ between fuzzy variable-basis $c$-spaces which are bijective,  fuzzy $c$-continuous and $c$-closed, forms a group. We can say that a way of seeing fuzzy variable-basis topology is studying invariants of the action of these groups aver the category $SET\times\mathcal{D}$.                           
\section{Idempotent and additive closure operators}
\begin{defi}
The closure operator $C=(c_{\text{\tiny{$X$}}})_{\text{\tiny{$X\in |SET|$}}}$ of definition (\ref{CO}) is called idempotent if the condition
\[
c_{\text{\tiny{$X$}}}\big(c_{\text{\tiny{$X$}}}(u)\big)= c_{\text{\tiny{$X$}}}(u)\quad \text{for all}\quad u\in L^X
\]
holds for every set $X$.
\end{defi}
\begin{prop}
Let $C=(c_{\text{\tiny{$Y$}}})_{\text{\tiny{$Y\in |SET|$}}}$ be an  idempotent closure operator. Then 
the initial closure operator $C=(c_{\text{\tiny{$X_f$}}})_{\text{\tiny{$X\in |SET|$}}}$ defined by 
\[
c_{\text{\tiny{$X_f$}}}=f_{\text{\tiny{$L$}}}^{\text{\tiny{$\leftarrow$}}} \circ c_{\text{\tiny{$Y$}}}\circ f_{\text{\tiny{$L$}}}^{\text{\tiny{$\rightarrow$}}} \qquad \text {for each function} \qquad f:X\rightarrow Y
\]
is also idempotent.
\end{prop}
\begin{proof}
Suppose that $C=(c_{\text{\tiny{$Y$}}})_{\text{\tiny{$Y\in |SET|$}}}$ is an  idempotent closure operator and let $f: X \rightarrow Y$ be a function. Then
\begin{align*}
c_{\text{\tiny{$X_f$}}}\circ c_{\text{\tiny{$X_f$}}}&=(f_{\text{\tiny{$L$}}}^{\text{\tiny{$\leftarrow$}}} \circ c_{\text{\tiny{$Y$}}}\circ f_{\text{\tiny{$L$}}}^{\text{\tiny{$\rightarrow$}}})\circ(f_{\text{\tiny{$L$}}}^{\text{\tiny{$\leftarrow$}}} \circ c_{\text{\tiny{$Y$}}}\circ f_{\text{\tiny{$L$}}}^{\text{\tiny{$\rightarrow$}}})\\
&=(f_{\text{\tiny{$L$}}}^{\text{\tiny{$\leftarrow$}}} \circ c_{\text{\tiny{$Y$}}})\circ( f_{\text{\tiny{$L$}}}^{\text{\tiny{$\rightarrow$}}}\circ f_{\text{\tiny{$L$}}}^{\text{\tiny{$\leftarrow$}}}) \circ( c_{\text{\tiny{$Y$}}}\circ f_{\text{\tiny{$L$}}}^{\text{\tiny{$\rightarrow$}}})\\
&\ls f_{\text{\tiny{$L$}}}^{\text{\tiny{$\leftarrow$}}} \circ (c_{\text{\tiny{$Y$}}}\circ c_{\text{\tiny{$Y$}}})\circ f_{\text{\tiny{$L$}}}^{\text{\tiny{$\rightarrow$}}}\\
&=f_{\text{\tiny{$L$}}}^{\text{\tiny{$\leftarrow$}}} \circ c_{\text{\tiny{$Y$}}}\circ f_{\text{\tiny{$L$}}}^{\text{\tiny{$\rightarrow$}}}\\
&=c_{\text{\tiny{$X_f$}}}.
\end{align*}
On the other hand, the monotonicity condition of closure operators implies that
\[c_{\text{\tiny{$X_f$}}}\ls c_{\text{\tiny{$X_f$}}}\circ c_{\text{\tiny{$X_f$}}}.
\]
\end{proof}
By using similar arguments, we can proof that 
\begin{prop}
Let 
$C=(c_{\text{\tiny{$XL$}}})_{\text{\tiny{$(X,L)$}} \in \left| SET\times\mathcal{D} \right|}$ be an idempotent closure operator. Then the initial closure operator
 $\hat{c}_{\text{\tiny{$XL$}}}:L^{\text{\tiny{$X$}}}\longrightarrow L^{\text{\tiny{$X$}}}$, defined by 
\[
 \hat{c}_{\text{\tiny{$XL$}}}=(f,\phi)^\leftarrow \circ c_{\text{\tiny{$YM$}}}\circ (f,\phi)^\rightarrow,
\]
for each morphism $(f,\phi):(X,L)\longrightarrow(Y,M)$ in $SET\times\mathcal{D}$ is also idempotent.
\end{prop}
\begin{defi}
The closure operator $C=(c_{\text{\tiny{$X$}}})_{\text{\tiny{$X\in |SET|$}}}$ of definition (\ref{CO}) is called 
\begin{enumerate}
\item Additive if the condition
\[
c_{\text{\tiny{$X$}}}(u\lor v)= c_{\text{\tiny{$X$}}}(u)\lor c_{\text{\tiny{$X$}}}(v)\quad \text{for all}\quad u,v\in L^X
\]
holds for every set $X$.
\item Fully additive if the condition
\[
c_{\text{\tiny{$X$}}}(\bigvee_{\lambda \in \Lambda}u_{\lambda})= \bigvee_{\lambda \in \Lambda}c_{\text{\tiny{$X$}}}(u_{\lambda})\quad \text{for all}\quad \{u_{\lambda}\mid \lambda \in \Lambda\}\subseteq L^X
\]
holds for every set $X$.
\end{enumerate}
\end{defi}
\begin{prop}
Let $C=(c_{\text{\tiny{$Y$}}})_{\text{\tiny{$Y\in |SET|$}}}$ be a fully additive closure operator. Then 
the initial closure operator $C=(c_{\text{\tiny{$X_f$}}})_{\text{\tiny{$X\in |SET|$}}}$ defined by 
\[
c_{\text{\tiny{$X_f$}}}=f_{\text{\tiny{$L$}}}^{\text{\tiny{$\leftarrow$}}} \circ c_{\text{\tiny{$Y$}}}\circ f_{\text{\tiny{$L$}}}^{\text{\tiny{$\rightarrow$}}} \qquad \text {for each function} \qquad f:X\rightarrow Y
\]
is also fully additive.
\end{prop}
\begin{proof}
Suppose that $C=(c_{\text{\tiny{$Y$}}})_{\text{\tiny{$Y\in |SET|$}}}$ is a fully additive closure operator and let $f: X \rightarrow Y$ be a function. Then,
$\text{for all}\quad \{u_{\lambda}\mid \lambda \in \Lambda\}\subseteq L^X$,
\begin{align*}
c_{\text{\tiny{$X_f$}}}(\bigvee_{\lambda \in \Lambda} u_{\lambda})&=f_{\text{\tiny{$L$}}}^{\text{\tiny{$\leftarrow$}}} \Big( c_{\text{\tiny{$Y$}}}\big( f_{\text{\tiny{$L$}}}^{\text{\tiny{$\rightarrow$}}}(\bigvee_{\lambda \in \Lambda}u_{\lambda})\big)\Big)   \\
&=\bigvee_{\lambda \in \Lambda}(f_{\text{\tiny{$L$}}}^{\text{\tiny{$\leftarrow$}}} \Big( c_{\text{\tiny{$Y$}}}\big( f_{\text{\tiny{$L$}}}^{\text{\tiny{$\rightarrow$}}}(u_{\lambda})\big)\Big)   \\
&=\bigvee_{\lambda \in \Lambda}c_{\text{\tiny{$X_f$}}}(u_{\lambda}).
\end{align*}
\end{proof}
By using similar arguments, we can proof that 
\begin{prop}
Let 
$C=(c_{\text{\tiny{$XL$}}})_{\text{\tiny{$(X,L)$}} \in \left| SET\times\mathcal{D} \right|}$ be a fully additive  closure operator. Then the initial closure operator   $\hat{c}_{\text{\tiny{$XL$}}}:L^{\text{\tiny{$X$}}}\longrightarrow L^{\text{\tiny{$X$}}}$, defined by 
\[
 \hat{c}_{\text{\tiny{$XL$}}}=(f,\phi)^\leftarrow \circ c_{\text{\tiny{$YM$}}}\circ (f,\phi)^\rightarrow,
\]
for each morphism $(f,\phi):(X,L)\longrightarrow(Y,M)$ in $SET\times\mathcal{D}$ is also fully additive.
\end{prop}
\section {Some examples of closure operators}
\begin{ex}
Let $X = \{x\}$ be a single point set and $L = [0, 1]$ be the usual unit interval. The  maps $c_{\text{\tiny{$n$}}}:L^X\rightarrow L^X$ defined by $c_{\text{\tiny{$n$}}}(t)=t^{\frac1n}$, for $n=1,2,\cdots$ are closure maps, from which just $c_{\text{\tiny{$1$}}}$ and $c_{\text{\tiny{$\infty$}}}=\lim\limits_{n\to \infty}t^{\frac1n}$ are idempotent.
\end{ex}
For a $GL-$monoid $\L$ and for an $\L$-topology $\tau\subseteq \L^X$, we define
\[
c_{\text{\tiny{$X$}}}(u)=\bigwedge_{v\in \tau}\{v\impl 0_X\mid  u\ls v\}.
\]
These  maps produce a closure operator of the category SET.
\begin{ex}
Let $X = \{x\}$ be a single point set and $\L$ be the set of all positive divisor of 12. It is clear that $\L$, ordered by 
\[ a\ls b \quad \text{ if and only if}\quad ''a\,\ \text{is divisor of}\,\ b'',
\]
is a  $GL-$monoid, where $a\otimes b= m.c.d\{a,b\}$. \\

If $\tau=\{1,2,12\}$ then
\[
c_{\text{\tiny{$X$}}}(1)=1; \quad c_{\text{\tiny{$X$}}}(2)=12;\quad c_{\text{\tiny{$X$}}}(3)=3; 
\]
\[\,\ c_{\text{\tiny{$X$}}}(4)=12;\quad c_{\text{\tiny{$X$}}}(6)=12;\quad c_{\text{\tiny{$X$}}}(12)=12
\]
produces an idempotent closure map.
\end{ex}
\begin{ex}
Let $X = \{x\}$ be a single point set and $\L$ be the set of all positive divisor of 36. $\L$, as in the previous example, ordered by 
\[ a\ls b \quad \text{ if and only if}\quad ''a\,\ \text{is divisor of}\,\ b'',
\]
is a  $GL-$monoid, where $a\otimes b= m.c.d\{a,b\}$. \\
The non-idempotent closure map $ c_{\text{\tiny{$X$}}}: \L\rightarrow \L$ defined by
\[
c_{\text{\tiny{$X$}}}(1)=1; \quad c_{\text{\tiny{$X$}}}(2)=2;\quad c_{\text{\tiny{$X$}}}(3)=6; 
\]
\[\,\ c_{\text{\tiny{$X$}}}(4)=4;\quad c_{\text{\tiny{$X$}}}(6)=18;\quad c_{\text{\tiny{$X$}}}(9)=9;
\]

\[\,\ c_{\text{\tiny{$X$}}}(12)= c_{\text{\tiny{$X$}}}(18)= c_{\text{\tiny{$X$}}}(36)=36
\]
produces an $\L$-topology 
\[
\tau =\{a\impl 1\mid c_{\text{\tiny{$X$}}}(a)= a\}=\{1,4,9,36 \}.
\]
\end{ex}
\begin{ex}
Let $\mu$ be a fuzzy subgroup of a group $G$.
 Acording to Theorem $1.3$ of \cite{YZ}, the fuzzy normalizer $N(\mu)$  of $\mu$ is a subgroup of $G$ and $\mu$ is a fuzzy normal subgroup of the group $N(\mu)$. Therefore 
\[
c(\mu)=\bigwedge\{ \eta\mid \eta \ls \mu\unlhd  G\}, \quad \text{for all  $G$ in the category GRP of all groups}
\] 
define the fuzzy normal closure operator of groups, which is idempotent.

\end{ex}

\end{document}